\documentclass[11pt]{amsart}
\usepackage[latin9]{inputenc}
\usepackage{amsthm}
\usepackage{amssymb}

\makeatletter
\numberwithin{equation}{section}
\numberwithin{figure}{section}
\theoremstyle{plain}
\newtheorem{thm}{\protect\theoremname}
  \theoremstyle{definition}
  \newtheorem{example}[thm]{\protect\examplename}
  \theoremstyle{plain}
  \newtheorem{cor}[thm]{\protect\corollaryname}
  \theoremstyle{plain}
  \newtheorem{lem}[thm]{\protect\lemmaname}

\title{}

\usepackage{amsfonts}\usepackage{enumerate}

\theoremstyle{plain}

\numberwithin{equation}{section}

\makeatother

  \providecommand{\corollaryname}{Corollary}
  \providecommand{\examplename}{Example}
  \providecommand{\lemmaname}{Lemma}
\providecommand{\theoremname}{Theorem}

\begin{document}

\title{Some properties of Stein-AG-groupoids adn Stein test}

\author{M. Rashad$^{a}$}

\address{Department of Mathematics, University of Malakand, Dir(L), Pakistan.}

\email{rashad@uom.edu.pk}

\author{I. Ahmad$^{a}$}

\email{iahmad@uom.edu.pk}

\address{a. Department of Mathematics University of Malakand, KPK, Pakistan.}

\author{M. Shah}

\address{Department of Mathematics, GPGC, Mardan, Pakistan. }

\email{shahmaths\_problem@hotmail.com}

\author{Amanullah}

\address{Department of Mathematics, University of Makand, Pakistan. }

\begin{abstract}
A groupoid that satisfying the left invertive law is called an AG-groupoid.
this concept is extended to introduce a Stein AG-groupoid. We prove
the existence  by providing some non-associative
examples. We also explore some basic and general properties of these
AG-groupoids and find their relations with other subclasses
of AG-groupoids.
\end{abstract}

\maketitle

\section{Introduction and priliminary}

A groupoid $(S,\cdot)$ or shortly $S$ is called\textit{ an AG-groupoid
}if it satisfied \textit{left invertive law: }$(ab)c=(cb)a$. This
structure is also known as\textit{ left almost semi-group }(LA-semigroup)
in \cite{key-9}, \textit{left invertive groupoid} in \cite{key-10},
while\textit{\ right modular groupoid} in \cite{key-11}. An AG-groupoid
$S$ is said to be an AG-Band if all of its elements are idempotent,
and is called AG-3-band if $a(aa)=(aa)a=a\mbox{\ensuremath{\forall}\ensuremath{a\ensuremath{\in}S.}}$
In {[}\ref{QM}{]} it is proved that $AG^{**}$-3-band is commutative
semigroup. 

In this paper we are going to invistigate new properties $AG^{*}$and
$AG^{**}.$ An AG-groupoid $S$ always satisfies the \textit{medial
law}: $(ab)(cd)=(ac)(bd)$ while an AG-groupoid $S$ with left identity
$e$ satisfies \textit{paramedial law}: $(ab)(cd)=(db)(ca)$. 

\section{Suggestion for detail studies}

An AG-groupoid $S$ is called \textit{transitively-commutative} if
$\forall\,a,b,c\in S,\,ab=ba,\,bc=cb\,\implies ac=ca$.

An AG-groupoid satisfying the identity $(ab)c=b(ac)$ is called $AG^{*}$-groupoid.
While if an AG-groupoid satisfying the identity $a(bc)=b(ac)$ is
called $AG^{**}$-groupoid.

In {[}\ref{MAI}, \ref{MIA}{]} some new classes of AG-groupoid has
been discovered that are $Bol^{*}$, $T^{1},T^{2},T^{3},T^{4}$ and
some other classes of AG-groupoids. Here we take the definitions of
$T^{1}$, $T^{3}$ and $T^{4}$ AG-groupoids to further invistigate
them.

An AG-groupoid $S$ is called \textit{$T^{1}$-AG-groupoid} if $\forall\,a,b,c,d\in S,\,\,\,\,ab=cd\,\implies ba=dc$.
An AG-groupoid $S$ is called \textit{left nuclear square }if $\forall\,a,b,c,\in S,\,\,\,\,a^{2}(bc)=(a^{2}b)c$.
Similarly $S$ is called right nuclear square if $\forall\,a,b,c,\in S,\,\,\,\,(ab)c^{2}=a(bc^{2})$
and middle nuclear square if $\forall\,a,b,c,\in S,\,\,\,\,(ab^{2})c=a(b^{2}c).$
An element $a$ of $S$ is called left cancellative if $ax=ay\implies x=y\,\,\forall x,y\in S.$
Similarly an element $a$ of $S$ is called right cancellative if
$xa=ya\implies x=y\,\forall x,y\in S.$ An element $a$ of $S$ is
called cancellative if it is both left and right cancellative. $S$
is called left cancellative (right cancellative, cancellative) if
every element of $S$ is left cancellative (right cancellative, cancellative).
A groupoid $S$ is called $Bol^{*}$-groupoid if it satisfiy the identity
$a(bc.d)=(ab.c)d$. And an AG-groupoid $S$ is called left alternative
if $\forall\,a,b\in S,\,\,\,\,aa.b=a.ab$ 

\section{properties of some subclasses of ag-groupoid}
\begin{thm}
\label{theorem-1}Every paramedial AG-groupoid $S$ is left nuclear
square.
\end{thm}
\begin{proof}
Let $S$ be a paramedial AG-groupoid, and let $a,b,c\in S$. Then
\begin{eqnarray*}
a^{2}(bc) & = & (aa)(bc)\\
 & = & ca.ba\,\,\qquad\mbox{by paramedial law}\\
 & = & cb.a^{2}\,\,\qquad\mbox{by medial law}\\
 & = & (a^{2}b)c\,\qquad\mbox{by left invertive law}\\
\implies\,a^{2}(bc) & = & (a^{2}b)c.
\end{eqnarray*}
Hence paramedial AG-groupoid is left nuclear square.
\end{proof}
Here is an example of left nuclear square AG groupoid that is not
paramedial AG groupoid.
\begin{example}
A left nuclear square AG-groupoid that is not paramedial AG-groupoid. 

\[
\begin{tabular}{l|llll}
 \ensuremath{\cdot}  &  \ensuremath{1}  &  \ensuremath{2}  &  \ensuremath{3}  &  \ensuremath{4}\\
\hline  \ensuremath{1}  &  \ensuremath{1}  &  \ensuremath{1}  &  \ensuremath{1}  &  \ensuremath{1}\\
 \ensuremath{2}  &  \ensuremath{1}  &  \ensuremath{1}  &  \ensuremath{1}  &  \ensuremath{3}\\
 \ensuremath{3}  &  \ensuremath{1}  &  \ensuremath{4}  &  \ensuremath{1}  &  \ensuremath{1}\\
 \ensuremath{4}  &  \ensuremath{1}  &  \ensuremath{1}  &  \mbox{2}  &  1 
\end{tabular}
\]

Clearly $(24)(43)\neq(34)(42).$
\end{example}
In {[}\ref{MAI}, Theorem 7{]} it is proved that every $T^{1}$-AG-groupoid
is $Bol^{*}$-AG-groupoid, and in {[}\ref{MIA}, Lemma 10{]} it is
proved that every $Bol^{*}$-AG-groupoid is paramedial AG-groupoid.
Now in this paper we prove that every $T^{1}$-AG-groupoid is paramedial
but the converse is not true.
\begin{thm}
Every $T^{1}$- AG-groupoid is paramedial-AG-groupoid.
\end{thm}
\begin{proof}
Proof: Let $S$ be a $T^{1}-$AG-groupoid, and let $a,b,c,d\in G$.
Then by definition of $T^{1}$-AG-groupoid.
\begin{eqnarray*}
ab=cd & \implies & ba=dc
\end{eqnarray*}
Now since,
\begin{eqnarray*}
ab.cd & = & ac.bd\,\,\qquad\mbox{by medial law}\\
\implies cd.ab & = & bd.ac\,\,\qquad\mbox{by definition of \ensuremath{T^{1}}-AG-groupoid}\\
 & = & (ac.d)b\qquad\mbox{by left invertive law}\\
 & = & (dc.a)b\qquad\mbox{by left invertive law}\\
 & = & ba.dc\,\,\qquad\mbox{by left invertive law}\\
\implies ab.cd & = & dc.ba\,\,\qquad\mbox{by definition of \ensuremath{T^{1}}-AG-groupoid}\\
 & = & db.ca\,\,\qquad\mbox{by medial law}\\
\implies ab.cd & = & db.ca
\end{eqnarray*}
Hence $S$ is paramedial-AG-groupoid.
\end{proof}
Here is an example of paramedial AG-groupoid that is not $T^{1}$-
AG-groupoid.
\begin{example}
Paramedial AG-groupoid of order 3 that is not $T^{1}-AG-$groupoid.
\[
\begin{tabular}{l|lll}
 \ensuremath{\cdot}  &  \ensuremath{1}  &  \ensuremath{2}  &  \ensuremath{3}\\
\hline  \ensuremath{1}  &  \ensuremath{1}  &  \ensuremath{1}  &  \ensuremath{1}\\
 \ensuremath{2}  &  \ensuremath{1}  &  \ensuremath{1}  &  \ensuremath{1}\\
 \ensuremath{3}  &  \ensuremath{2}  &  \ensuremath{2}  &  \ensuremath{2} 
\end{tabular}
\]
\end{example}
Clearly
\[
(1*2)=(2*3)\implies(2*1)\neq(3*2)
\]

Next, we prove that every $AG^{*}$-groupoid is $Bol^{*}$-AG-groupoid
and hence by {[}\ref{MIA}, Lemma 10{]} $Bol^{*}-$AG-groupoid is
paramedial, so one can conclude that $AG^{*}$-groupoid is paramedial.
Also by Theorem (\ref{theorem-1}) $AG^{*}$-groupoid is left nuclear
square.
\begin{thm}
Every $AG^{*}$-groupoid is $Bol^{*}$-AG-groupoid.
\end{thm}
\begin{proof}
Let $S$ be an $AG^{*}$-groupoid, and let $a,b,c,d\,\in G$. Then
by definition of $AG^{*}$-groupoid
\begin{eqnarray*}
(ab)c & = & b(ac)
\end{eqnarray*}
Now
\begin{eqnarray*}
(ab.c)d & = & dc.ab\,\,\qquad\mbox{(by left invertive law})\\
 & = & da.cb\,\,\qquad\mbox{(by medial law)}\\
 & = & a(d.cb)\qquad\mbox{(by definition of \ensuremath{AG^{*}\mbox{)}}}\\
 & = & a(cd.b)\qquad\mbox{(by definition of \ensuremath{AG^{*}\mbox{)}}}\\
 & = & a(bd.c)\qquad\mbox{(by left invertive law)}\\
 & = & a(d.bc)\qquad\mbox{(by definition of \ensuremath{AG^{*}\mbox{)}}}\\
 & = & da.bc\,\,\qquad\mbox{(by definition of \ensuremath{AG^{*}\mbox{)}}}\\
 & = & (bc.a)d\qquad\mbox{(by left invertive law)}\\
 & = & a(bc.d)\qquad\mbox{(by definition of \ensuremath{AG^{*}\mbox{)}}}\\
\implies(ab.c)d & = & a(bc.d).
\end{eqnarray*}
Hence $AG^{*}$-groupoid is $Bol^{*}$-AG-groupoid.
\end{proof}
Here is a counterexample that $Bol^{*}$-AG-groupoid is not $AG^{*}$-groupoid.
\begin{example}
$Bol^{*}$-AG-groupoid of order 3 that is not $AG^{*}$groupoid.
\[
\begin{tabular}{l|lll}
 \ensuremath{\cdot}  &  \ensuremath{1}  &  \ensuremath{2}  &  \ensuremath{3}\\
\hline  \ensuremath{1}  &  \ensuremath{1}  &  \ensuremath{1}  &  \ensuremath{1}\\
 \ensuremath{2}  &  \ensuremath{1}  &  \ensuremath{1}  &  \ensuremath{1}\\
 \ensuremath{3}  &  \ensuremath{1}  &  \ensuremath{2}  &  \ensuremath{1} 
\end{tabular}
\]

Clearly
\[
(3*3)*2\neq3*(3*2).
\]
\end{example}
\begin{cor}
Every $AG^{*}$-groupoid is paramedial AG-groupoid.
\end{cor}
\begin{proof}
Let $S$ be an $AG^{*}$-groupoid and let $a,b,c,d\,\in G$. Then
by definition of $AG^{*}$-groupoid
\begin{eqnarray*}
(ab)c & = & b(ac)
\end{eqnarray*}
Now,
\begin{eqnarray*}
ab.cd & = & (cd.b)a\qquad\mbox{(by left invertive law)}\\
 & = & (d.cb)a\qquad\mbox{(by definition of \ensuremath{AG^{*}})}\\
 & = & cb.da\,\,\qquad\mbox{(by definition of \ensuremath{AG^{*}})}\\
 & = & (da.b)c\qquad\mbox{(by left invertive law)}\\
 & = & (a.db)c\qquad\mbox{(by definition of \ensuremath{AG^{*}})}\\
 & = & (c.db)a\qquad\mbox{(by left invertive law)}\\
 & = & db.ca\,\,\qquad\mbox{(by definition of \ensuremath{AG^{*}})}\\
\implies ab.cd & = & db.ca.
\end{eqnarray*}
Thus $S$ is Paramedial AG-groupoid.
\end{proof}

\begin{cor}
Every $AG^{*}$-groupoid is left nuclear square AG-groupoid.
\end{cor}
\begin{proof}
Let $S$ be an $AG^{*}$-groupoid, and let $a,b,c\in S.$ Then 
\begin{eqnarray*}
(a^{2}b)c & = & cb.\,\,a^{2}\qquad\mbox{by left invertive law}\\
 & = & b(c.aa)\qquad\mbox{by definition of \ensuremath{AG^{*}}}\\
 & = & b(ac.a)\qquad\mbox{by definition of \ensuremath{AG^{*}}}\\
 & = & (ac.b)a\qquad\mbox{by definition of \ensuremath{AG^{*}}}\\
 & = & (bc.a)a\qquad\mbox{by left invertive law}\\
 & = & a^{2}(bc)\,\qquad\mbox{by left invertive law}\\
\implies(a^{2}b)c & = & a^{2}(bc).
\end{eqnarray*}
Hence, every $AG^{*}$-groupoid is left nuclear. 
\end{proof}
\begin{thm}
Every right cancellative $AG^{*}$-groupoid $S$ is transitively commutative
AG-groupoid.
\end{thm}
\begin{proof}
Let $S$ be an $AG^{*}$-groupoid, and let $a,b,c\in S$ such that
\[
ab=ba\,\,\,\,\,\mbox{and}\,\,\,\,\,\,\,\,bc=cb.
\]
Then consider
\begin{eqnarray*}
(ac)b & = & c(ab)\qquad\mbox{(by definition of \ensuremath{AG^{*}}-groupoid})\\
 & = & c(ba)\qquad\mbox{(as}\,\mbox{\ensuremath{ab=ba})}\\
 & = & (bc)a\qquad\mbox{(by definition of \ensuremath{AG^{*}}-groupoid}\\
 & = & (cb)a\qquad\mbox{\mbox{(as}\,\mbox{\ensuremath{bc=cb})}}\\
 & = & (ab)c\qquad\mbox{(by left invertive law)}\\
 & = & (ba)c\qquad\mbox{(as\,}\mbox{\ensuremath{ab=ba})}\\
 & = & (ca)b\qquad\mbox{(by left invertive law)}\\
\implies(ac)b & = & (ca)b\\
ac & = & ca\,\,\,\,\qquad\mbox{(by right cancellativity)}
\end{eqnarray*}
Hence right cancellative $AG^{*}$-groupoid is transitively commutative-AG-groupoid.
\end{proof}
\begin{thm}
Every $T^{1}$-AG-3-band is $AG^{*}$-groupoid.
\end{thm}
Proof: Let $S$ be a $T^{1}$-AG-groupoid, and let $a,b,c,d\in G$.
Then by definition of $T^{1}$-AG-groupoid
\begin{proof}
\begin{eqnarray*}
ab=cd & \implies & ba=dc
\end{eqnarray*}
Now since,
\begin{eqnarray*}
(ab)c & = & (cb)a\,\,\,\,\,\,\,\,\qquad\mbox{(by left invertive law)}\\
 & = & (cb)((aa)a)\qquad\mbox{(by definition of AG-3-band)}\\
 & = & c(aa).ba\,\,\,\,\qquad\mbox{(by medial law)}\\
\implies c(ab) & = & ba.c(aa)\,\,\,\,\qquad\mbox{(by definition of \ensuremath{T^{1}}-AG-groupoid)}\\
 & = & ((c(aa))a)b\,\qquad\mbox{(by left invertive law)}\\
\implies(ab)c & = & b((c(aa))a)\,\qquad\mbox{(by definition of \ensuremath{T^{1}}-AG-groupoid})\\
 & = & b((a(aa)c)\,\,\qquad\mbox{(by left invertive law)}\\
 & = & b(ac)\,\,\,\,\,\,\,\,\,\qquad(\mbox{by definition of AG-3-band)}\\
\implies(ab)c & = & b(ac)
\end{eqnarray*}
Thus $S$ is $AG^{*}$-groupoid.
\end{proof}
\begin{lem}
Every $AG^{*}$-groupoid is left alternative.
\end{lem}
\begin{proof}
Let $S$ be an $AG^{*}$-groupoid, and let $a,b,c\in S.$ Then by
definition of $AG^{*}$-groupoid, we have
\begin{eqnarray}
(ab)c & = & b(ac)\label{eq:a}
\end{eqnarray}
Now replacing $b$ by $a$ in (\ref{eq:a}), we have
\begin{eqnarray*}
(aa)b & = & a(ab)
\end{eqnarray*}
Hence $S$ is left alternative.
\end{proof}
\begin{lem}
Every $AG^{*}$-groupoid is middle nuclear square AG-groupoid.
\end{lem}
\begin{proof}
Let $S$ be an $AG^{*}$-groupoid, and let $a,b,c\in S.$ Then by
definition of $AG^{*}$-groupoid, we have
\begin{eqnarray*}
(ab)c & = & b(ac)
\end{eqnarray*}
Then
\begin{eqnarray*}
(ab^{2})c & = & b^{2}.ac\,\,\,\,\qquad\mbox{(by definition of \ensuremath{AG^{*}})}\\
 & = & bb.ac\\
 & = & ba.bc\,\,\,\,\qquad\mbox{(by medial law)}\\
 & = & a(b.bc)\,\,\qquad\mbox{(by definition of \ensuremath{AG^{*}})}\\
 & = & a(bb.c)\,\,\qquad\mbox{(by definition of \ensuremath{AG^{*}})}\\
 & = & a(b^{2}c)\\
\implies(ab^{2})c & = & a(b^{2}c).
\end{eqnarray*}
Hence $S$ is middle nuclear square AG-groupoid.
\end{proof}
\begin{thm}
Every $AG^{*}$-groupoid is right nuclear square AG-groupoid.
\end{thm}
\begin{proof}
Let $S$ be an $AG^{*}$-groupoid, and let $a,b,c\in S.$ Then by
definition of $AG^{*}$-groupoid, we have
\begin{eqnarray*}
(ab)c & = & b(ac)
\end{eqnarray*}
Now consider,
\begin{eqnarray*}
a(bc^{2}) & = & ba.c^{2}\,\,\,\,\,\qquad\mbox{(by definition of \ensuremath{AG^{*}})}\\
 & = & bc.ac\,\,\,\,\,\qquad\mbox{(by medial law)}\\
 & = & (ac.c)b\,\,\,\qquad\mbox{(by left invertive law)}\\
 & = & (cc.a)b\,\,\,\qquad\mbox{(by left invertive law)}\\
 & = & a(cc.b)\,\,\,\qquad\mbox{(by definition of \ensuremath{AG^{*}})}\\
 & = & a(c.cb)\,\,\,\qquad\mbox{(by definition of \ensuremath{AG^{*}})}\\
 & = & ca.cb\,\,\,\,\,\qquad\mbox{(by definition of \ensuremath{AG^{*}})}\\
 & = & c^{2}.ab\,\,\,\,\,\qquad\mbox{(by medial law)}\\
 & = & (ab.c)c\,\,\,\qquad\mbox{(by left invertive law)}\\
 & = & (b.ac)c\,\,\,\qquad\mbox{(by definition of \ensuremath{AG^{*}})}\\
 & = & ac.bc\,\,\,\,\,\qquad\mbox{(by definition of \ensuremath{AG^{*}})}\\
 & = & ab.c^{2}\,\,\,\,\,\qquad\mbox{(by medial law)}
\end{eqnarray*}

Hence $S$ is right nuclear square AG-groupoid.
\end{proof}
In {[}\ref{MIA}, Lemma 9{]} it is proved that every $AG^{**}$-groupoid
$S$ is $Bol^{*}$-groupoid and by {[}\ref{MIA}, Lemma 10{]} every
$Bol^{*}$-groupoid is paramedial AG- groupoid. Hence immediatly we
have following result . 
\begin{lem}
Every $AG^{**}$groupoid is paramedial-AG-groupoid.
\end{lem}
Prrof: Let $S$ be an $AG^{**}$-groupoid, and let $a,b,c,d\in G.$
Then
\begin{eqnarray*}
ab.cd & = & c(ab.d)\qquad\mbox{(}\mbox{by \ensuremath{AG^{**}}}-\mbox{groupoid)}\\
 & = & c(db.a)\qquad\mbox{\mbox{(}by medial law)}\\
 & = & db.ca\,\,\qquad\mbox{(}\mbox{by \ensuremath{AG^{**}}}-\mbox{groupoid)}
\end{eqnarray*}
Hence $S$ is  paramedial-AG-groupoid.

\end{document}